\documentclass[11pt]{article}
\usepackage{amsfonts}
\usepackage{nicefrac}
\usepackage{amsmath,amscd}

\begin{document}

\newtheorem{theorem}{Theorem}[section]
\newtheorem{lemma}[theorem]{Lemma}
\newtheorem{proposition}[theorem]{Proposition}
\newtheorem{corollary}[theorem]{Corollary}
\newtheorem{definition}[theorem]{Definition}
\newtheorem{example}[theorem]{Example}

\newenvironment{proof}[1][Proof]{\begin{trivlist}
\item[\hskip \labelsep {\bfseries #1}]}{\end{trivlist}}
\newenvironment{remark}[1][Remark]{\begin{trivlist}
\item[\hskip \labelsep {\bfseries #1}]}{\end{trivlist}}

\newcommand{\qed}{\nobreak \ifvmode \relax \else
      \ifdim\lastskip<1.5em \hskip-\lastskip
      \hskip1.5em plus0em minus0.5em \fi \nobreak
      \vrule height0.75em width0.5em depth0.25em\fi}

\title{Positivstellens\"atze for semi-algebraic sets in real closed valued fields}
\author{Noa Lavi \\ \footnotesize{Department of Mathematics} \\ \footnotesize{The Hebrew University of Jerusalem} \\ \footnotesize {Givat Ram, Jerusalem 91904, Israel} \\                 \footnotesize{noa.lavi@mail.huji.ac.il} }
\maketitle

\begin{center}
    \scriptsize The purpose of this paper is to give a characterization for polynomials and rational functions which admit only non-negative values on definable sets in real closed valued fields. That is, generalizing the relative positivstellens\"atze for sets defined also by valuation terms. For this, we use model theoretic tools, together with existence of canonical valuations.
\end{center}
\section {Introduction}
A ``Nichtnegativstellensatz'' in real algebraic geometry is a theorem 
which gives an algebraic characterization of those polynomials admitting
 only non-negative values on a given set. The original Nichtnegativstellensatz 
is the solution by Artin in \cite{art}, Theorem 45, to Hilbert's seventeenth problem: a polynomial taking only non-negative values can be written as a sum of squares 
of rational functions. Later on A.Robinson \cite{rob} generalized the theorem to any real closed field using the model completeness of real closed fields \cite{trsk}, a model theoretic property which means that any formula is equivalent to an existential formula in the language of rings. For proof see \cite{PLNM} Theorem 5.1 p.49, and Theorem 5.7 p.54. \\
In real algebra, valuations come very naturally into the picture. A real closed field $K$ admits a canonical valuation which is non-trivial if and only if the field is non-archimedeam, that is, not embeddable into $\mathbb{R}$. The  canonical valuation ring is defined to be the convex hull of $\mathbb{Z}$ inside $K$, and its maximal ideal is exactly the set of the infinitesimal elements known from nonstandard analysis. \\
An ordered valued field is an ordered field with a convex valuation ring. In this paper we will consider a real closed field with some convex valuation and its field of multi-variable rational functions. For some real closed valued field $K$, one could consider semi-algebraic sets defined by valuation inequalities (in addition to ordinary polynomial inequalities), which we also refer to as ``valuative semi-algebraic sets'' in this paper. That is, \[S_{\bar q}=\Big\{\bar x \in K^n | \nu(\bar q(\bar x))\ge 0 \Big\}\] where $\bar q \subset K(\bar X)$. The valuation ring is an example of such a set, as it is defined by the formula $\nu(\bar X)\ge 0$. In \cite{Dick}, Dickmann proved what may be paraphrased as the following: 
\begin{theorem} \cite{Dick} Theorem $2$, p.132. \\
Let $K$ be a real closed valued field, let ${\cal O}_K$ be its valuation ring and ${\cal M}$ its set of infinitesimal elements. Then for every $p \in {\cal O}_K[\bar X]$, $p$ is non-negative on ${\cal O}_K$ if and only if \[p=\sum \frac{f_i^2(1+m_it_i)}{g_i^2(1+n_iw_i)} \] where $f_i,g_i,t_i,w_i \in {\cal O}_K[\bar X]$ and $m_i,n_i \in {\cal M}$.
\end{theorem} 
 In this paper we generalize the above result by proving a relative Positivstellensatz for valuative semi-algebraic sets. This could be seen as the analogue of the work done in \cite{LY}, and was, as a matter a fact, the reason to work on the problem at first place.
\begin{theorem} \label{mainpos}
Let $K$ be a real closed valued field, and let $f_1,...f_k, h_1,...,h_l \in K(\bar X), g_1,...,g_{k'}, u_1,...,u_{l'} \in K[\bar X]$. We define the valuative semi-algebraic set \[
S_{\bar f, \bar h, \bar g, \bar u}=\Big\{ \bar x \in K^n | \nu(\bar f(\bar x)) \ge 0, \nu(\bar h(\bar x))>0, g(\bar x) \ge 0, u(\bar x) > 0 \Big\}.
\]
Then for every $p \in K[\bar X]$, $p$ is non-negative on $S_{\bar f, \bar h, \bar g, \bar u}$ if and only if \[p=\sum \frac{a_i^2c_i(1+m_it_i)(1+s_i)}{b_i^2d_i(1+n_iw_i)(1+r_i)} \] where $a_i,b_i \in K[\bar X]$, $t_i,w_i$ are in the ${\cal O}_K$-algebra generated by $\bar f$ and $ \bar h$, and $c_i,d_i$ are in the multiplicative monoid generated by $\bar g, \bar u$, and $s_i,r_i$ are in the multiplicative monoid generated by ${\cal O}_K$ and $\bar h$, and $m_i,n_i \in {\cal M}$.
\end{theorem}
We also prove a strict Positivstellensatz for such sets.
\begin{theorem} \label{strictpos}
Let $K$ be a real closed valued field, and let $f_1,...f_k, h_1,...,h_l \in K(\bar X), g_1,...,g_{k'}, u_1,...,u_{l'} \in K[\bar X]$. For every $p \in K[\bar X]$, we have $p$ strictly positive on $S_{\bar f, \bar h, \bar g, \bar u}$ if and only if \[p=\sum \frac{1+a_i^2c_i(1+m_it_i)(1+s_i)}{b_i^2d_i(1+n_iw_i)(1+r_i)} \] where $a_i,b_i \in K[\bar X]$, $t_i,w_i$ are in the ${\cal O}_K$-algebra generated by $\bar f$ and $ \bar h$, and $c_i,d_i$ are in the multiplicative monoid generated by $\bar g, \bar u$, and $s_i,r_i$ are in the multiplicative monoid generated by ${\cal O}_K$ and $\bar h$, and $m_i,n_i \in {\cal M}$. 
\end{theorem}
We finish by characterizing the ordered fields that have such Positivstellensatz, and prove that those are exactly the ordered fields which are dense in their real closure. 
\section { Preliminaries}
We begin by defining the class of fields which we shall work on in this paper.
\begin{definition}
An \emph{ordered valued field} (or $OVF$) $ \langle K,\nu,\le \rangle$ is a valued field with an order satisfying: for every $x,y \in K$ if $0<x\le y$ then $\nu (y)\le \nu(x) $.
\end{definition}
An \emph{$OVF$-valuation} will be a valuation $\nu$ on a field $K$ such that there exist some ordering $\le_K$ on $K$ such that $\langle K, \nu, \le_K \rangle \models OVF$.
\begin{definition}
A \emph{real closed valued field} (or $RCVF$) is an order valued field which is also real closed.
\end{definition}
A main tool that we use is the relative canonical valuation of an ordered field. 
\begin{definition} \label{canon}
Let $K$ be a field and $P$ be the set of all the positive elements of $K$ according to some order $\le_P$ on $K$. The \emph{canonical valuation} induced by $P$, denoted by ${\cal O}_P$, is the ring of all elements $x \in K$ such that there exists $n \in \mathbb{N}$ such that $-n \le_P x \le_P n$. Let $F\subset K$ a subfield with ${\cal O}_F$ a valuation ring which makes $F$ together with the induced order from $\le_P$ an ordered valued field. Then the \emph{canonical valuation relative to ${\cal O}_F$} is the ring of all elements $x \in K$ such that there exist $a,b \in {\cal{O}}_F$ such that $a \le_P x \le_P b$, i.e, the convex hull of ${\cal O}_F$ in $K$. We denote it by $Val_{{\cal O}_F}(P)$.
\end{definition}
We recall also the definition of the cone generated by a subset of a field.
\begin{definition} \cite{PD}, Definition 1.1.6, p.10. \\ \label{cone}
Let $L$ be some commutative ring with a unit, $T\subseteq L$ some subset.  The \emph{cone} generated by $T$ is the minimal set containing $T\cup (L^{\times})^2 $ which is closed under addition and multiplication. We denote it here as \emph{$Cone(T)$}. For such a cone $P$ we denote $P \cap -P$ by $suppP$. When $P \cup -P = L$ and $suppP$ is a prime ideal, we say that $P$ is a \emph{positive cone}. The \emph{real spectrum} of $L$, denoted by $Sper{L}$ is the set of all positive cones of $L$.
\end{definition}
In order to prove our theorem, we use elementary properties of ordered and real closed valued fields. In \cite{CD}, Theorem 4B, P.218, G. Cherlin and M. Dickmann proved that $RCVF$ is the model companion of $OVF$. This means that if ${\cal L} = \langle L, \tilde\nu, \le_L \rangle \models OVF$ is an extension of the real closed valued field ${\cal K} = \langle K, \nu, \le_K \rangle$, then for every formula $\phi(\bar {X}) $ with parameters from $K$, if there exists $\bar v \in L^n$ such that ${\cal{L}}  \models \phi(\bar v)$ then there exists $\bar b \in K^n$ such that ${\cal{K}} \models \phi(\bar b) $. This plays the role of the model completeness of $Th(\mathbb{R})$ in the proof of Hilbert's 17th problem.
\begin{definition}
Let $L_{ovring}$ denote the language of $L_{ring}$ enriched by symbols for order and valuation. We call a quantifier free formula in $L_{ovring}$ a \emph{generalized semi-algebraic formula} if it is a conjunction of order and valuation inequalities. That is, in addition to the ordinary order inequalities, we have also conditions of the form $\nu(p(\bar X) \ge \nu(q(\bar X)) $ and $\nu(p(\bar X)) > \nu(q(\bar X)) $.
\end{definition}

\section{Positivstellens\"atze for valuative semi algebraic sets over real closed valued fields}
By \cite{art}, we already know that if $K$ is a real closed field, then a polynomial $p \in K[\bar X]$ is non-negative on $K^n$ if and only if there exists $f_1,...,f_k \in K(\bar X)$ such that $p=\sum_{i=1}^k {f_i}^2$. Suppose $K$ is a real closed valued field. It is easy to see that $1+mp({\bar x})$ is non-negative on for every $m \in \cal{M}$ (=the maximal ideal of ${\cal O}_K$), $\bar x \in {\cal O}_K^n $ and $p(\bar X) \in {\cal O}_K[\bar X]$. In \cite{Dick}, Dickmann has shown that these sum up to all the positive elements on ${\cal O}_K$ among ${\cal O}_K[\bar X]$. In order to generalize the above result, we start by proving the following lemma: 
\begin{lemma} \label {FOP}
Let $K$ be an ordered valued field, $f \in K(\bar X)$ and $\bar b \in K^n$. If $f$ is defined on $\bar b$ we have: 

 $f(\bar b) \in {\cal O}_K$ if and only if $1+mf(\bar b) > 0$ for every $m \in {\cal M}_K$.
\end{lemma} 
\begin{proof}
One direction is obvious. For the other direction, suppose $f(\bar b) \notin {\cal O}_K$. Without loss of generality we assume $f(\bar b) > 0$. Since $\nu(\frac{1}{f(\bar b)})>0$, taking $m = -\frac{1}{f(\bar b)} \in {\cal M}_K$ we have $1+mf(\bar b)=0$, contradiction. $\qed$
\end{proof}
\begin{remark}
We note that the same thing works if we require $f(\bar b) \in {\cal M}_K$ and $m \in {\cal O}_K$ instead of $f(\bar b) \in {\cal O}_K$ and $m \in {\cal M}_K$.
\end{remark}
We now quote the general abstract Positivstellensatz theorem which we use to prove the results of this paper.
\begin{theorem} \cite{PD}, Theorem 4.2.1, p.86. \\ \label{AP}
Let $A$ be some commutative ring with a unit. For subsets $F,G,H \subset A$, let $I(F)$ denote the ideal generated by $F$ and $\prod G$ the multiplicative monoid of $A$ generated by $G$. Then, the following are equivalent: \\
(i) There is no $P \in SperA$, for which the following conditions hold simultaneously: $F \subseteq suppP$, $G \subseteq A - suppP$, $ H \subseteq P$. \\
(ii) There are $b \in I(F)$, $c \in \prod G$, and $t \in Cone(H)$ such that $c^2+t=b$. 
\end{theorem}
Now we generalize the transfer principle at Theorem $4.2.8$ at \cite{PD} so it will fit also for real closed valued fields.
\begin{theorem} \label{transfer}
Let $K$ be a real closed valued field and let $\phi(\bar X)$ a generalized semi-algebraic formula in the language $L_{ovring}$ with parameters from $K$. There exists some $\bar b \in K^n$ such that ${\cal K} \models \phi(\bar b)$ if and only if there exists some $P \in Sper(K[\bar X])$ such that \[\langle \nicefrac{K[\bar X]}{suppP},\nu_P, \le_{ \nicefrac{P}{suppP} } \rangle \models \Large \phi(\nicefrac{\bar X}{suppP}) \]
where $\nu_P$ is the canonical valuation $Val_{{\cal O}_K}(\le_{ \nicefrac{P}{suppP} }) $.
\end{theorem}
\begin{proof}
For one direction, let $\bar b \in K^n$ be such that ${\cal K} \models \phi(\bar b) $. Let us define $P$ exactly as in the proof of transfer principle in Theorem $4.2.8$ in \cite{PD}, that is,  $P := \{a  \in K[\bar X] | a(\bar b) \ge 0 \}$. For atomic sub-formula of $\phi$ of the form ``$v \ge 0$'' or ``$v > 0$'', it is clear. For every atomic sub-formulas of $\phi$ of the form ``$\nu(v) \ge \nu(w)$'' we obtain by Lemma \ref{FOP} $\nu(p(\bar b)) \ge \nu(q(\bar b))$ if and only if for every $m \in {\cal M} $ we have $q(q+mp)(\bar b) > 0 $, hence if and only if $q(q+mp) \in P $ and $q(q+mp) \notin suppP $. Same thing for any atomic sub-formula of the form ``$\nu(v)> \nu(w)$''. Again, by Lemma \ref{FOP} we get \[\langle \nicefrac{K[\bar X]}{suppP},\nu_P, \le_{ \nicefrac{P}{suppP} } \rangle \models \Large \phi( \nicefrac {\bar X}{suppP}). \] \\
For the other direction, we move to the real closure of $\langle \nicefrac{K[\bar X]}{suppP}, \tilde{\le}_{ \nicefrac{P}{suppP} } \rangle$, that we shall denote by $\tilde{K}$. Let $\tilde{\nu_P}$ and $\tilde{\le}_{ \nicefrac{P}{suppP} }$ be the corresponding extensions of $\nu_P$ and $\le_{ \nicefrac{P}{suppP} }$ to $\tilde{K}$. Then we have \[
\langle \tilde{K},\tilde{\nu_P}, \tilde{\le}_{ \nicefrac{P}{suppP} } \rangle \models \Large \phi(\nicefrac{\bar X}{suppP}).
\] 
Hence, from model-completeness of $RCVF$ \cite{CD}, there exists some $\bar b \in K^n$ such that $\langle K, \nu, \le \rangle \models \phi(\bar b)$. $\qed$
\end{proof}
We are now ready to prove the main result of the paper. 

\begin{proof} [Proof of Theorem \ref{mainpos}.]
For one direction, a straightforward computation using Lemma \ref{FOP} shows that any such polynomial is, indeed, non-negative on $S_{\bar f,\bar h, \bar g, \bar u}$. For the other direction, let \[H=\Big\{1+m_if_i | 1\le i \le k, m_i \in {\cal M} \Big\} \cup \Big\{1+o_ih_i | 1\le i \le l, o_i \in {\cal O}_K\Big\} \cup Q  \cup \Big\{-p\Big\}\] where $Q=\Big\{g_1,\ldots,g_{k'},u_1,\ldots,u_{l'}\Big\} $,  
$G=\Big\{u_1,\ldots, u_{l'} \Big\} \cup \Big\{p\Big\}$ and $F=\Big\{0\Big\}$. Then according to Theorem \ref{AP} the following are equivalent: \\ (i) There is no $P \in Sper(K[\bar X])$ such that $H \subseteq P$ and $p \notin suppP$ \\ (ii) There exists $c \in \Pi G$ and $t \in Cone(H)$ such that $c^2+t=0$. \\
 Hence, (i) is equivalent to the existence of \[t_1, t_2  \in Cone(\Big\{1+m_if_i | 1\le i \le k, m_i \in {\cal M} \Big\} \cup \Big\{1+o_ih_i | 1\le i \le l, o_i \in {\cal O}_K\Big\} \cup Q),\] and $c \in \prod{G}$ such that $c^2+t_1-pt_2=0$. That is, $p=\frac{c^{2}+t_1}{t_2}. $ As $c^{2}$ is a square and hence belongs to $Cone(H)$, $p$ has, indeed, the required form. If $p$ is not of this form then there exists $P \in Sper(K[\bar X])$ such that \[\Big\{1+m_if_i | 1\le i \le k, m_i \in {\cal M} \Big\} \cup \Big\{1+o_ih_i | 1\le i \le l, o_i \in {\cal O}_K\Big\} \cup Q \subseteq P\] as well as $-p \in P$ and $p \notin suppP$. Hence  \[
\langle \tilde{K},\tilde{\nu_P}, \tilde{\le}_{ \nicefrac{P}{suppP} } \rangle \models \nu(\bar f(\nicefrac{\bar X}{suppP})) \ge 0 \wedge \bar \nu ( h(\nicefrac{\bar X}{suppP})) > 0  \wedge \bar g((\nicefrac{\bar X}{suppP})) \ge 0 \] \[ \wedge \bar u((\nicefrac{\bar X}{suppP})) \ge 0 \wedge \bar u((\nicefrac{\bar X}{suppP})) \neq 0 \wedge p(\nicefrac{\bar X}{suppP}) <0.
\] 
Hence, by Theorem \ref{transfer}, $p$ cannot be
non-negative on $S_{\bar f, \bar h, \bar g, \bar u}$ if it doesn't have required form. $\qed$
\end{proof}
We may also deduce a strict Positivstellensatz in the very same way.
\begin{proof} [Proof of Theorem \ref{strictpos}.]
We let $H$ and $F$ be the same as in the proof of Theorem \ref{mainpos}, and $G= \Big\{u_1,\ldots, u_{l'} \Big\} \cup \Big \{1\Big\}$. Then we get $p=\frac{c^2+t_1}{t_2}$, where $c \in \prod_G$. As $c^2 \in Cone(H)$ we get $p=\frac{1+t_1}{t_2}$. The proof continues as in the proof of Theorem \ref{mainpos}. $\qed$
\end{proof}

\subsection {Ordered valued fields with Positivstellensatz}
In \cite{MCK}, Theorem 2, p. 223, McKenna has characterized the fields admitting the characterization in Hilbert 17th problem, denoted by ``Hilbert property''. He showed that those fields are exactly those ordered fields which are dense in their real closure.  Following him, we call the characterization of Theorem \ref{mainpos} on $S_{\bar f, \bar h, \bar u}$ (that is, open sets), the ``Relative Hilbert Property'' (RHP), and show that an ordered valued field has this property if and only if it is dense in its real closure.  
\begin{proposition}
Let $\langle K, \nu, \le \rangle \models OVF$. Then $K$ has RHP if and only if $K$ is dense in its real closure.
\end{proposition}
\begin{proof}
Suppose $K$ is not dense in its real closure, denoted by $\tilde{K}$. In \cite{MCK} it is proved that there exists $p \in K[ x]$ such that $p$ is positive definite on $K$ but not in $\tilde{K}$. Hence, there exists $a_1,a_2 \in K$ such that $p$ changes sign in $[a_1,a_2]$ inside $\tilde{K}$ but not inside $K$. Let $B_{\gamma}$ be a closed ball containing $[a_1,a_2]$, where $\gamma \in \Gamma_K$. As $B_{\gamma}$ is defined by the formula $"\nu(x) \ge \nu(c) "$ where $c \in K$ is such that $\nu(c)=\gamma$, then $p$ cannot be of the form in Theorem \ref{mainpos}, as this would imply that $p$ is positive definite on $B_{\gamma}$ in $\tilde{K}$. Hence, $K$ doesn't have RHP. \\
Suppose $K$ is dense in its real closure $\tilde{K}$. Let $S \subseteq K^n$ be a set of the form in Theorem \ref{mainpos} and assume $p \in K[\bar X]$ is not of the form in Theorem \ref{mainpos}. Applying Theorem \ref{AP} with suitably chosen sets $F,G H \subset K[\bar X]$, there exists some $P \in Sper(K[\bar X])$ such that $\le_{\nicefrac{P}{suppP}}$ makes $p$ negative and all the polynomials appearing in the definition of $S$ positive. The order $\le_{\nicefrac{P}{suppP}}$ extends, of course, to $\tilde{K}[\bar X]$. Due to density of $K$ in $\tilde{K}$, for every $m \in {\cal M}_{\tilde{K}}$ we have $1+mf_i \tilde{\ge}_{\nicefrac{P}{suppP}} 0$ and for every $o \in {\cal O}_{\tilde{K}}$ we have $1+or_i \tilde{\ge}_{\nicefrac{P}{suppP}} 0$. Hence, by Theorem \ref{transfer}, $p$ is not positive definite on $S_{\bar f, \bar h, \bar u}(\tilde{K})$. Due to density, it is  not positive definite on $S_{\bar f, \bar h, \bar u}({K})$ either. $\qed$
\end{proof}


\addcontentsline{toc}{section}{References}
\begin{thebibliography}{9}

\bibitem{art} E. Artin, 
\"Uber die Zerlegung definiter Funktionen in Quadrate,
 {\it Abh. Sem. Univ. Hamburg} {\bf 5} (1927), 100-115.


\bibitem{CD}  G. Cherlin and M. Dickmann,
Real closed rings.  II. Model Theory,
 {\it Ann. Pure Appl. Logic} {\bf 25} (1983), no. 3, 213-231.

\bibitem{Dick} M. Dickmann,
On polynomials over real closed rings,
{\it Model theory of Algebra and Arithmetic, Proceeding Karpacz 1979 (L. Pacholski, J. Wiezejewski and A.Wilkie, eds.)}, Springer Lecture Notes in Mathematics 843 (1980), pp. 117-135

 \bibitem{LY} N. Lavi,
A Ganzstellensatz for open semi-algebraic sets in real
closed valued fields, {\it arXiv:1101.2116v1 }

\bibitem{MCK} K. McKenna,
New facts about Hilbert's seventeenth problem. {\it In: Model theory and Algebra , A memorial tribute to Abraham Robinson (D. Saracino and V. Weispfenning, eds.) . } Springer lecture notes in mathematics {\bf 498} (1975), pp. 220 - 230.

\bibitem{PLNM} A. Prestel, 
{\it Lectures on formally real fields}, Springer Lecture Notes in Mathematics 1093 (1984).

\bibitem{PD} A. Prestel, and C.N. Delzell,
 Positive polynomials: from Hilbert's 17th problem to real algebra. Springer, 2001.

\bibitem{rob} A. Robinson,
On ordered fields and definite functions.
{\it Mathematische Annalen 130 (1956),257-271.}

\bibitem{trsk} A. Tarski,
A Decision Method for Elementary Algebra and Geometry.
{\it University of California Press, Berkeley 1948.}


\end{thebibliography}
\end{document}